\documentclass[12pt,a4paper]{article}
\usepackage{amssymb}
\usepackage{amsfonts}

\usepackage{amsmath}
\usepackage{amsthm}
\usepackage{enumerate}

\RequirePackage[OT1]{fontenc}
\RequirePackage{amsthm,amsmath, xcolor}
\RequirePackage[numbers]{natbib}
\RequirePackage[colorlinks,citecolor=blue,urlcolor=blue]{hyperref}
\usepackage{amscd,amsfonts,amssymb,latexsym,array,hhline,graphicx}
\usepackage{float}

\newtheorem{theorem}{Theorem}[section]
\newtheorem{proposition}[theorem]{Proposition}

\newtheorem{remark}{Remark}[section]
\newtheorem{corollary}[theorem]{Corollary}

\newcommand\cA{{\cal A}}

\newcommand\cG{{\cal G}}

\newcommand\cL{{\cal L}}
\newcommand\cB{{\cal B}}

\newcommand\cD{{\cal D}}
\newcommand\cQ{{\cal Q}}

\newcommand\cR{{\cal R}}

\newcommand\ve{\varepsilon}

\newcommand\Er{\mbox{Err}}

\def\bbr{{\mathbb R}}

\def\text#1{\hbox{#1}}

\def\E{{\bf E}}

\def\B{{\bf B}}

\def\C{{\bf C}}
\def\D{{\bf D}}

\def\G{{\bf G}}
\def\L{{\bf L}}
\def\U{{\bf U}}

\def\c{{\bf c}}

\def\r{{\bf r}}

\newcommand\Tr{\mbox{Tr}}
\def\Chi{{\bf 1}}

\def\d{\mathrm{d}}
\def\build #1_#2{\mathrel{\mathop{\kern 0pt #1}\limits_{#2}}}
\newcommand\tr{\mbox{tr}}

\newcommand{\wh}{\widehat}
\newcommand{\wt}{\widetilde}
\newcommand{\zs}[1]{{\mathchoice{#1}{#1}{\lower.25ex\hbox{$\scriptstyle#1$}}
{\lower0.25ex\hbox{$\scriptscriptstyle#1$}}}}

\numberwithin{equation}{section}

\begin{document}
\title{
Oracle inequalities for the stochastic differential equations
\thanks{
This work was done under
 the Ministry of Education and Science of the Russian Federation in the framework of the research project no. 2.3208.2017/4.6,
 by RFBR Grant 16-01-00121 A
and by the
 partial financial support
of the RSF grant number 14-49-00079 (National Research University "MPEI"
14 Krasnokazarmennaya, 111250 Moscow, Russia).
}}
\author{Pchelintsev E.A.,
\thanks{
Department of Information Technologies and Business Analytics,
 Tomsk State University,
Lenin str. 36,
 634050 Tomsk, Russia,
 e-mail: evgen-pch@yandex.ru
}
 \and
 Pergamenshchikov S.M.\thanks{
 Laboratory of Mathematics LMRS,
  University of Rouen, France,
International Laboratory SSP \& QF,
 National Research Tomsk State University
and
National Research University "MPEI"
14 Krasnokazarmennaya, 111250 Moscow, Russia,
 e-mail: Serge.Pergamenchtchikov@univ-rouen.fr}
}
 \date{}
\maketitle

\begin{abstract}
This paper is a survey of recent results 
on the adaptive robust non parametric
methods for the continuous time regression model  with the semi - martingale  noises with jumps.
The noises are modeled by the L\'evy processes, the 
Ornstein -- Uhlenbeck processes and semi-Markov processes.
We represent the general model selection method and the sharp oracle inequalities methods 
which provide the robust efficient estimation in the adaptive setting.
Moreover, we present the recent results on the improved model selection methods for the nonparametric estimation problems.  
\end{abstract}

\\
\par
{\bf AMS (2010) Subject Classification : primary 62G08; secondary 62G05}

\bibliographystyle{plain}

\newpage

\section{Introduction}\label{sec:In}

This paper is a survey on the adaptive 
non parametric estimation
methods 
for the general semi-martingale regression model in continuous time  
defined as
\begin{equation}\label{sec:In.1}
 \d\,y_\zs{t} = S(t)\d\,t + \,
\d \xi_\zs{t}\,,\quad
 0\le t \le n\,,
\end{equation}
where $S(\cdot)$ is an unknown $1$ - periodic
 function,
$(\xi_\zs{t})_\zs{0\le t\le n}$  is an
unobservable noise defined by  semimartingale  with the values in the
Skorokhod space $\cD[0,n]$
 such that, for any
function $f$ from $\cL_\zs{2}[0,n]$, the stochastic integral
\begin{equation*}\label{sec:In.2}
I_\zs{n}(f)=\int^n_\zs{0}f(s)\d \xi_\zs{s}
 \end{equation*} is
well defined and has the following properties
\begin{equation}\label{sec:In.3}
\E_\zs{Q} I_\zs{n}(f)=0 \quad\mbox{and}\quad \E_\zs{Q}
I^2_\zs{n}(f)\le \kappa_\zs{Q} \int^{n}_\zs{0}\,f^{2}(s)\d s
\,.
\end{equation}
 We use $\E_\zs{Q}$ for the expectation with respect to the distribution
 $Q$ in $\cD[0,n]$ of the process $(\xi_\zs{t})_\zs{0\le t\le n}$,
 which is assumed to belong to some
 probability family $\cQ_\zs{n}$ and $\kappa_\zs{Q}$ is
some positive constant depending on the distribution $Q$.  
The problem consists to estimate the function $S$ on the observations $(y_\zs{t})_\zs{0\le t\le n}$.
Note that if $(\xi_\zs{t})_\zs{0\le t\le n}$ is a  brownian motion, then we obtain the well known "signal+white noise" model
which is very popular in statistical radio-physics (see, for example, \cite{IbragimovKhasminskii1981, Kutoyants1977, Kutoyants1984, Pinsker1981}).
 In this paper we assume that in addition
 to the intrinsic noise in the radio-electronic system, approximated usually by the Gaussian white or color noise,
 the useful signal $S$ is distorted by the impulse flow described by the processes with jumps.
The cause of the appearance of a pulse stream in the radio-electronic systems can be, for example, either external unintended (atmospheric) or intentional impulse noise
 and
the errors in the demodulation and the channel decoding for the binary information symbols. Note that,
for the first time the impulse noises
for the detection signal problems
 have been introduced
 on the basis of the compound Poisson processes
was introduced  by Kassam in
\cite{Kassam1988}.
 Later, such processes was used
  in
  \cite{KonevPergamenshchikov2012, KonevPergamenshchikov2015, KPP2014, Pchelintsev2013}
for the
parametric and
nonparametric signal estimation problems.
 However,  the compound Poisson process can  describe only   the large impulses  influence of fixed  single frequency.
 Taking into account  that in the telecommunication systems,  the impulses are without limitations on frequencies
 one needs to extend the framework of the observation model by making use the  L\'evy
 processes  \eqref{sec:In.1+1} which is a particular case of the general semimartinagale regression model
 introduced in
\cite{KonevPergamenshchikov2009a}.
Generally, we
 consider  nonparametric estimation problems for the function $S$ from
$\cL_\zs{2}$ under the condition that 
the distribution of the noise $(\xi_\zs{t})_\zs{0\le t\le n}$ is unknown. We know only that this distribution belongs to 
 some distribution family $\cQ_\zs{n}$. In this case
we use the robust estimation approach 
proposed  in \cite{GaltchoukPergamenshchikov2006, KonevPergamenshchikov2012, KonevPergamenshchikov2015}
for the nonparametric estimation. According to this approach we have to construct an estimator $\wh{S}_\zs{n}$
 (any function of $(y_\zs{t})_\zs{0\le t\le n}$) for $S$ to minimize the robust risk defined as
\begin{equation}\label{sec:risks}
\cR^{*}_\zs{n}(\wh{S}_\zs{n},S)=\sup_\zs{Q\in\cQ_\zs{n}}\,
\cR_\zs{Q}(\wh{S}_\zs{n},S)\,,
\end{equation}
where  $\cR_\zs{Q}(\cdot,\cdot)$ is the usual quadratic risk of the form
\begin{equation}\label{sec:risks_00}
\cR_\zs{Q}(\wh{S}_\zs{n},S):=
\E_\zs{Q}\,\|\wh{S}_\zs{n}-S\|^2
\quad\mbox{and}\quad
\Vert S\Vert^{2}=\int^{1}_\zs{0}\,S^{2}(t)\d t
\,.
\end{equation}  
It is clear that if we don't know the distribution
of the observation one  needs to find an estimator which will be optimal for all possible observation distributions. 
Moreover in this paper we consider the estimation problem in the adaptive setting,
i.e. when the regularity of $S$ is unknown. 
To this end
 we use the adaptive method based on the model selection approach. 
The interest to such statistical procedures is explained by the fact that
they provide
adaptive solutions for the nonparametric estimation through
oracle inequalities which give the  non-asymptotic upper
bound for the quadratic risk including
the minimal risk over chosen  family of estimators.
It should be noted that for the first time the model selection methods
were proposed by
 Akaike \cite{Akaike1974} and Mallows \cite{Mallows1973}
for parametric models. Then,
these methods had been developed
for the nonparametric estimation
 and
the oracle inequalities for the quadratic risks was obtained  by
Barron, Birg\'e and Massart \cite{BarronBirgeMassart1999},
Massart \cite{Massart2005},
by Fourdrinier and Pergamenshchikov \cite{FourdrinierPergamenshchikov2007}
for the regression models in discrete time
and
\cite{KonevPergamenshchikov2010} in continuous time.
Unfortunately, the oracle inequalities obtained in these papers
 can not
provide the efficient estimation in the adaptive setting, since
the upper bounds in these inequalities
have some fixed coefficients in the main terms which are more than one.
To obtain  the efficiency property for estimation
procedures  one has to obtain
 the sharp oracle inequalities, i.e.
 in which
  the factor at the principal term on the right-hand side of the inequality
 is close to unity.
The first result on sharp inequalities is most likely due to Kneip \cite{Kneip1994}
who studied a Gaussian regression model in the discrete time.
It will be observed that the
derivation of oracle inequalities usually rests upon
the fact that the initial model, by applying the Fourier transformation,
can be reduced to the Gaussian independent observations.
However, such  transformation is possible only for
Gaussian models with independent homogeneous observations or for
 inhomogeneous ones with  known correlation characteristics.
For the general non Gaussian observations one needs to use the approach proposed by
 Galtchouk and Pergamenshchikov
 \cite{GaltchoukPergamenshchikov2009a, GaltchoukPergamenshchikov2009b}  for the heteroscedastic
 regression models in discrete time and developed then by
Konev and Pergamenshchikov in
\cite{KonevPergamenshchikov2009a, KonevPergamenshchikov2009b, KonevPergamenshchikov2012, KonevPergamenshchikov2015}
for semimartingale models in  continuous time.
In general the model selection is an
adaptive rule $\wh{\alpha}$ which choses an estimator $S^{*}=\wh{S}_\zs{\wh{\alpha}}$ from 
an estimate family 
$(\wh{S}_\zs{\alpha})_\zs{\alpha\in\cA}$. The goal of  this selection is to prove the following
nonasymptotic 
oracle inequality: for any sufficient small $\delta >0$ and any observation duration $n\ge 1$
\begin{equation}\label{sec:In.3+5}
\cR( S^{*},S)\,
\le\,(1+\delta)\,\min_\zs{\alpha\in\cA}\,\cR(\wh{S}_\zs{\alpha},S)+\delta^{-1}
\cB_\zs{n}\,,
\end{equation}
where the rest term
$\cB_\zs{n}$  
 is sufficiently small with respect to the minimax convergence rate.
Such oracle inequalities are called {\sl sharp}, since the coefficient in the main term 
$1+\delta$
is close to one for sufficiently small $\delta>0$. 
Moreover, in this paper we represent the new results on the improved estimation methods for the
nonparametric models \eqref{sec:In.1}.
Usually, the model selection procedures are based on the least squares estimators.  But in
\cite{PchelintsevPchelintsev_VPergamenshchikov2017}
 it is propose
 to use the improved least square estimators which enable  to improve considerably the non asymptotic estimation accuracy.
At the first time such idea  was proposed in \cite{FourdrinierPergamenshchikov2007}
for the regression in discrete time and in
\cite{KonevPergamenshchikov2010}
for the Gaussian regression model in continuous time. In
\cite{PchelintsevPchelintsev_VPergamenshchikov2017}
 these methods 
are developed
for the non - Gaussion
regression models
in continuous time.
It should be noted that generally for the conditionally Gaussian regression models  one can not use the well known improved estimators proposed in
\cite{JamesStein1961, FourdrinierStrawderman1996} for Gaussian or spherically symmetric observations. To apply the improved estimation methods to the
 non Gaussian regression models in continuous time
 one needs to modify the well known James--Stein procedure in the way  proposed  in  \cite{Pchelintsev2013, KPP2014}.
 For the improved model selection procedures the oracle inequality \eqref{sec:In.3+5}
  is shown also. We note that this
inequality allows us to provide the asymptotic efficiency without knowing the regularity of the function being estimated. 
The efficacy property for a nonparametric estimate $S^{*}$ means
$$
\lim_\zs{n\to\infty}\,
\upsilon_\zs{n}\,\sup_\zs{S\in W^k_\zs{r}}\,
\cR(S^{*},S)\,
=
\lim_\zs{n\to\infty}\,
\upsilon_\zs{n}
\,
\inf_\zs{\wh{S}}
\sup_\zs{S\in W^k_\zs{r}}\,
\cR(\wh{S},S)\,
= 
l_\zs{*}(\r)\,,
$$
where
\begin{equation}\label{sec:Ae.3}
l_\zs{*}(\r)\,=\,((1+2k)\r)^{1/(2k+1)}\,
\left(\frac{k}{\pi (k+1)}\right)^{2k/(2k+1)}
\,,
\end{equation}
$\upsilon_\zs{n}$ is a normalizing coefficient (convergence rate),
 $W^{k}_\zs{r}$ is the Sobolev ball of a radius  $r>0$ and the regularity  $k\ge 1$. The limit \eqref{sec:Ae.3}
 is called the Pinsker constant 
 which is calculated by Pinsker in \cite{Pinsker1981}.

The rest of the paper is organized as follows. In the next section \ref{Examples}, we describe
the L\'evy, Ornstein--Uhlenbeck and semi-Markov processes as the examples of a
semimartingale impusle noise in the model \eqref{sec:In.1}. In Section \ref{sec:Mo} we construct the model selection
procedure based on the least square estimators and show the sharp oracle inequalities.
In Section  \ref{sec:Imp} we give the improved least squares estimators and
we study the improvement effect for the semimartingale model \eqref{sec:In.1}. In Section \ref{sec:IMo} we construct the improved model selection
procedure and show the sharp oracle inequalities. The asymptotic efficiency is studed in Section \ref{sec:Ae}.

\section{Examples}\label{Examples}

\subsection{L\'evy model}

First we consider the model \eqref{sec:In.1} with the L\'evy noise process, i.e. we assume that
 the noise process $(\xi_\zs{t})_\zs{0\le t\le n}$
is defined as
\begin{equation}\label{sec:In.1+1}
\xi_\zs{t}=\varrho_\zs{1} w_\zs{t} + \varrho_\zs{2} z_\zs{t}
\quad\mbox{and}\quad
z_\zs{t}=x*(\mu-\wt{\mu})_\zs{t}
\,,
\end{equation}
where $\varrho_\zs{1}$ and $\varrho_\zs{2}$ are  some unknown constants,
$(w_\zs{t})_\zs{t\ge\,0}$ is a standard brownian motion,
 $\mu(\d s\,\d x)$ is a jump measure with  deterministic
compensator $\wt{\mu}(\d s\,\d x)=\d s\Pi(\d x)$,
$\Pi(\cdot)$ is a L\'evy measure, i.e.  some positive measure on $\bbr_\zs{*}=\bbr\setminus \{0\}$,
(see, for example,
\cite{JacodShiryaev2002, ContTankov2004} for details) such that
\begin{equation*}\label{sec:Ex.1-00_mPi}
\Pi(x^{2})=1
\quad\mbox{and}\quad
\Pi(x^{6})
\,<\,\infty\,.
\end{equation*}
Here we use the notation $\Pi(\vert x\vert^{m})=\int_\zs{\bbr_\zs{*}}\,\vert z\vert^{m}\,\Pi(\d z)$. Note that the L\'evy measure
 $\Pi(\bbr_\zs{*})$ could be equal to $+\infty$.
 One can check directly that for the process \eqref{sec:In.1+1}
the condition \eqref{sec:In.3} holds with  $\kappa_\zs{Q}=\sigma_\zs{Q}=\varrho^{2}_\zs{1}+\varrho^{2}_\zs{2}$.
We assume that the nuisance parameters
    $\varrho_\zs{1}$
and $\varrho_\zs{2}$ of the process $(\xi_\zs{t})_\zs{0\le t\le n}$
 satisfy the conditions
\begin{equation}\label{sec:Ex.01-1}
0< \underline{\varrho}\le \varrho^{2}_\zs{1}
\quad\mbox{and}\quad
\sigma_\zs{Q}\,
\le
\varsigma^{*}
\,,
\end{equation}
where
 the bounds
 $\underline{\varrho}$
and
$\varsigma^{*}$
are functions of $n$, i.e.
 $\underline{\varrho}= \underline{\varrho}_\zs{n}$
and
$\varsigma^{*}=\varsigma^{*}_\zs{n}$
 such that for any $\check{\delta}>0$
\begin{equation}\label{sec:Ex.01-2}
\liminf_\zs{n\to\infty}\,n^{\check{\delta}}\,
\underline{\varrho}_\zs{n}
\,
>0
\quad\mbox{and}\quad
\lim_\zs{n\to\infty}\,n^{-\check{\delta}}\,\varsigma^{*}_\zs{n}
=0
\,.
\end{equation}

For this example $\cQ_\zs{n}$ is  the family
of all distributions of process \eqref{sec:In.1} -- \eqref{sec:In.1+1} on the Skorokhod space
$\D[0,n]$ satisfying the conditions \eqref{sec:Ex.01-1} -- \eqref{sec:Ex.01-2}.

The models \eqref{sec:In.1}
with the L\'evy's type noise
are used in different applied problems
(see \cite{Bertoin1996}, for details). Such  models naturally arise
in the nonparametric functional statistics problems
(see, for example, \cite{FerratyVieu2006}).

\subsection{Ornstein -- Uhlenbeck model}

Now we consider the noise process $(\xi_\zs{t})_\zs{t\ge 0}$ in
\eqref{sec:In.1} difened by a non-Gaussian Ornstein-Uhlenbeck
process with the L\'evy subordinator. Let the noise
process in \eqref{sec:In.1} obey the equation
\begin{equation}\label{sec:Ex.1}
\d\xi_\zs{t} = a\xi_\zs{t}\d t+\d u_\zs{t}\,,\quad \xi_\zs{0}=0\,,
\end{equation}
where
$
u_\zs{t}=
\varrho_\zs{1}\, w_\zs{t}+\varrho_\zs{2}\,z_\zs{t}$ and
the process
$z_\zs{t}$
is defied in \eqref{sec:In.1+1}.
Here $a\le 0$,
$\varrho_\zs{1}$ and $\varrho_\zs{2}$ are unknown parameters.
 We assume that the parameters  $\varrho_\zs{1}$
and $\varrho_\zs{2}$ satisfy the conditions
\eqref{sec:Ex.01-1} and the parameter
\begin{equation}\label{sec:Ex.3+}
-a_\zs{max}\le  a\le 0\,,
\end{equation}
where the bound $a_\zs{max}>0$ is the function of $n$, i.e. $a_\zs{max}=a_\zs{max}(n)$,
such that for any $\check{\delta}>0$
\begin{equation}\label{sec:Ex.3-1}
\lim_\zs{n\to\infty}\,\frac{a_\zs{max}}{n^{\check{\delta}}}=0
\,.
\end{equation}
In this case  $\cQ_\zs{n}$ is the family
of all distributions of process  \eqref{sec:Ex.1} on the Skorokhod space
$\D[0,n]$ satisfying the conditions \eqref{sec:Ex.01-1},  \eqref{sec:Ex.3+} and \eqref{sec:Ex.3-1}.
Note also that the processes
\eqref{sec:In.1+1} and
\eqref{sec:Ex.1} are $\cG$ - conditionally  Gaussian square integrated semimartingales,  where
 $\cG=\sigma\{z_\zs{t}\,,\,t\ge 0\}$.

Such processes are used in
the financial Black-Scholes type markets with  jumps (see, for
example, \cite{BarndorffNielsenShephard2001, DelongKluppelberg2008} and the references therein). 
Note also that in the case when $\varrho_\zs{2}=0$ for the parametric estimation problem such models are considered in
\cite{HopfnerKutoyants2009, HopfnerKutoyants2010, KonevPergamenshchikov2003}.

\subsection{Semi -- Markov model}

In 
\cite{BarbuBeltaifPergamenshchikov2017a, BarbuBeltaifPergamenshchikov2017b}
it is introduced
the regression model   \eqref{sec:In.1}  in which       
 the noise process describes by the equation
\begin{equation}\label{sec:Ex.5}
 \xi_\zs{t} = \varrho_\zs{1} L_\zs{t} + \varrho_\zs{2} X_\zs{t}\,,
\end{equation}
and  the L\'evy  process $L_\zs{t}$ is defined as 
 \begin{equation*}\label{sec:Mcs.1}
L_\zs{t}=\check{\varrho}\,w_\zs{t}
+
\sqrt{1-\check{\varrho}^{2}}\,z_\zs{t}\,,
\end{equation*}
 where $0\le \check{\varrho}\le 1$ is an unknown constant.
Moreover, we assume that the pure jump process $(X_\zs{t})_\zs{t\ge\, 0}$ in
 \eqref{sec:Ex.5}
 is a semi-Markov process with the following form
 \begin{equation}\label{sec:Ex.2}
 X_\zs{t} = \sum_\zs{i=1}^{N_\zs{t}} Y_\zs{i},
\end{equation}
where $(Y_\zs{i})_\zs{i\ge\, 1}$ is an i.i.d. sequence of random variables with
\begin{equation*}\label{sec:Ex.3+?++}
\E\,Y_\zs{i}=0\,,\quad 
\E\,Y^2_\zs{i}=1
\quad\mbox{and}\quad
\E\,Y^4_\zs{i}<\infty
\,.
\end{equation*}
Here $N_\zs{t}$ is a general  counting process (see, for example, \cite{Mikosch2004})
defined as
\begin{equation*}\label{sec:Ex.4}
N_\zs{t} = \sum_\zs{k=1}^{\infty} \Chi_\zs{\{T_\zs{k} \le t\}}
\quad\mbox{and}\quad
T_\zs{k}=\sum_\zs{l=1}^k\, \tau_\zs{l}\,,
\end{equation*}
where $(\tau_\zs{l})_\zs{l\ge\,1}$ is an i.i.d. sequence of positive integrated
 random variables with distribution $\eta$ and mean $\check{\tau}=\E\,\tau_\zs{1}>0$. We assume that the processes 
$(N_\zs{t})_\zs{t\ge 0}$ and  $(Y_\zs{i})_{i\ge\, 1}$ are independent between them and are also independent of $(L_\zs{t})_\zs{t\ge 0}$.
Here, the family  $\cQ_\zs{n}$ is defined by 
of all distributions of process  \eqref{sec:Ex.5} on the Skorokhod space
$\D[0,n]$ with the parameters $\varrho_\zs{1}$ and 
$\varrho_\zs{2}$
 satisfying the conditions \eqref{sec:Ex.01-1}  and $0\le \check{\varrho}\le 1$.

Note that the process $(X_\zs{t})_\zs{t\ge\, 0}$ is a special case of a semi-Markov process (see, e.g., \cite{BarbuLimnios2008} and \cite{LimniosOprisan}).
 It should be noted that if $\tau_\zs{j}$ are exponential random variables, then $(N_\zs{t})_\zs{t\ge 0}$
is a Poisson process and, in this case, $(\xi_\zs{t})_\zs{t\ge 0}$ is a L\'evy process. 
But, in the general case when the process \eqref{sec:Ex.2} is not a L\'evy process, this process has a memory and cannot be treated
in the framework of semi-martingales with independent increments.  In this case, we need to 
develop new tools based on renewal theory arguments from  \cite{Goldie1991}.

It should be noted  that for $\check{\varrho}>0$  the process \eqref{sec:Ex.5}
 is $\cG$ - conditionally  Gaussian also. In this case
 $\cG=\sigma\{z_\zs{t}\,,\,X_\zs{t}\,,\,t\ge 0\}$.

\section{Model selection}\label{sec:Mo}

\noindent 
Let $(\phi_\zs{j})_\zs{j\ge\, 1}$ be an orthonormal uniformly bounded basis in $\cL_\zs{2}[0,1]$, i.e. 
for some  constant $\phi_\zs{*}\ge 1$, which may be depend on $n$,
\begin{equation*}\label{sec:In.3-00}
\sup_\zs{0\le j\le n}\,\sup_\zs{0\le t\le 1}\vert\phi_\zs{j}(t)\vert\,
\le\,
\phi_\zs{*}
<\infty\,. 
\end{equation*}
For example, we can take 
 the trigonometric basis    defined as $\Tr_\zs{1}\equiv 1$ and, for $j\ge 2,$
\begin{equation}\label{sec:In.5}
 \Tr_\zs{j}(x)= \sqrt 2
\left\{
\begin{array}{c}
\cos(2\pi[j/2] x)\, \quad\mbox{for even}\quad j;\\[4mm]
\sin(2\pi[j/2] x)\quad\mbox{for odd}\quad j,
\end{array}
\right.
\end{equation}
where $[x]$ denotes the integer part of $x$.

To estimate the function $S$ we use here the model selection procedure for continuous time regression models from
\cite{KonevPergamenshchikov2012} based on the Fourrier expansion. We recall that  
for any function $S$ from $\cL_\zs{2}[0,1]$ we can write 
\begin{equation*}\label{sec:In.5+++Fourrier}
S(t)=\sum^{\infty}_\zs{j=1}\,\theta_\zs{j}\,\phi_\zs{j}(t)
\quad\mbox{and}\quad
\theta_\zs{j}= (S,\phi_\zs{j}) = \int_\zs{0}^{1} S(t) \phi_\zs{j}(t)\d t
\,.
\end{equation*}
So, to estimate the function $S$ it suffices to estimate the coefficients $\theta_\zs{j}$ and to replace them in this representation by their estimators.
Using the fact that the function $S$ and $\phi_\zs{j}$ are $1$ - periodic we can write that
$$
\theta_\zs{j}=\frac{1}{n} \int_\zs{0}^{n}\, \phi_\zs{j}(t)\,S(t) \d t
\,.
$$
If we replace here the differential $S(t)\d t$ by the stochastic observed differential $\d y_\zs{t}$
then we obtain the natural estimate for $\theta_\zs{j}$ on the time interval $[0,n]$ 
\begin{equation*}\label{sec:In.7}
\wh{\theta}_\zs{j,n}= \frac{1}{n} \int_\zs{0}^{n}  \phi_\zs{j}(t) \d\,y_\zs{t}\,,
\end{equation*}
which
can be represented, in view of the model \eqref{sec:In.1}, as
\begin{equation*}\label{sec:In.8}
\wh{\theta}_\zs{j,n}= \theta_\zs{j} + \frac{1}{\sqrt n}\xi_\zs{j,n}
\quad\mbox{and}\quad 
\xi_\zs{j,n}
= \frac{1}{\sqrt n} I_\zs{n}(\phi_\zs{j})\,.
\end{equation*}
We need to impose some stability conditions for the noise Fourier transform sequence 
$(\xi_\zs{j,n})_\zs{1\le j\le n}$. To this end we set  for some stability noise intensity parameter $\sigma_\zs{Q}> 0$
the following function
\begin{equation}
\label{L_1_Q}
\L_\zs{1,n}(Q)=
 \sup_\zs{x\in [-1,1]^{n}}\,
 \left|
\sum^{n}_\zs{j=1}\,x_\zs{j}\,
\left(
\E_\zs{Q}\,\xi^{2}_\zs{j,n}
-
\sigma_\zs{Q}
\right)
\right|
\,.
\end{equation}
In \cite{KonevPergamenshchikov2012}
the  parameter $\sigma_\zs{Q}$ is called proxy variance .

\noindent $\C_\zs{1})$ {\it There exists a variance proxy
$\sigma_\zs{Q}> 0$ 
such that for any $\epsilon>0$
$$
\lim_\zs{n\to\infty}\frac{\L_\zs{1,n}(Q)}{n^{\epsilon}}
=0\,.
$$
}

\noindent Moreover, we set 
\begin{equation*}
\label{L_2_Q}
\L_\zs{2,n}(Q)=
 \sup_\zs{|x|\le 1}
\E_\zs{Q}\, \left(
\sum^{n}_\zs{j=1}\,x_\zs{j}\,
(\xi^2_\zs{j,n}-\E_\zs{Q} \xi^2_\zs{j,n}) 
\right)^2\,.
\end{equation*}
\noindent $\C_\zs{2})$ {\it Assume that for any $\epsilon>0$
$$
\lim_\zs{n\to\infty}\frac{\L_\zs{2,n}(Q)}{n^{\epsilon}}
=0\,.
$$
}

Now (see, for example, \cite{IbragimovKhasminskii1981}) we can estimate the function $S$ by the projection estimators, i.e. 
\begin{equation}\label{sec:In.7++pre}
\wh{S}_\zs{m}(t)=\sum^{m}_\zs{j=1}\,\wh{\theta}_\zs{j,n}\,\phi_\zs{j}(t)\,,\quad 0\le t\le 1\,,
\end{equation}
for some number $m\to\infty$ as $n\to\infty$. It should be noted that Pinsker in \cite{Pinsker1981} shows that the projection estimators of the form
\eqref{sec:In.7++pre}
 are not efficient. For obtaining efficient estimation one needs to use weighted least square estimators defined as 
\begin{equation}\label{sec:Mo.1+??+}
\wh{S}_\lambda (t) = \sum_\zs{j=1}^{n} \lambda(j) \wh{\theta}_\zs{j,n} \phi_\zs{j}(t)\,,
\end{equation}
where the coefficients $\lambda=(\lambda(j))_\zs{1\le j\le n}$ belong to some finite set $\Lambda$ from $[0,1]^n$.
 As it is shown in  \cite{Pinsker1981}, in order to obtain efficient estimators, the coefficients $\lambda(j)$ in \eqref{sec:Mo.1+??+} need to be chosen depending on 
the regularity of the unknown function $S$. Since we consider the adaptive case, i.e. we assume that the regularity of the function $S$ is unknown, then we chose the weight coefficients on the basis of  the model selection procedure proposed in \cite{KonevPergamenshchikov2012}
   for the general semi-martingale regression model in continuous time. 
 To the end, first we set
 \begin{equation*}\label{sec:Mo.2++???}
\nu=\#(\Lambda)
\quad\mbox{and}\quad
\vert\Lambda\vert_\zs{*}=1+ \max_\zs{\lambda\in\Lambda}\,
\sum^{n}_\zs{j=1}\lambda(j)
\,,
\end{equation*}
where $\#(\Lambda)$ is the cardinal number of  $\Lambda$. 
Now, to choose a weight sequence $\lambda$ in the set $\Lambda$ we use the empirical  quadratic risk, defined as
$$
\Er_n(\lambda) = \parallel \wh{S}_\lambda-S\parallel^2,
$$
which in our case is equal to
\begin{equation*}\label{sec:Mo.3}
\Er_n(\lambda) = \sum_\zs{j=1}^{n} \lambda^2(j) \wh{\theta}^2_\zs{j,n} -2 \sum_\zs{j=1}^{n} \lambda(j) \wh{\theta}_\zs{j,n}\theta_\zs{j}+ \sum_\zs{j=1}^{\infty} \theta^2_\zs{j}.
\end{equation*}
Since the Fourier coefficients $(\theta_\zs{j})_\zs{j\ge\,1}$ are unknown, we replace
the terms $\wh{\theta}_\zs{j,n}\theta_\zs{j,n}$ by  
\begin{equation*}\label{sec:Mo.4+?/}
\wt{\theta}_\zs{j,n} = \wh{\theta}^2_\zs{j,n} - \frac{  \wh{\sigma}_\zs{n}}{n}\,,
\end{equation*}
where $\wh{\sigma}_\zs{n}$ is an estimate for the variance proxy $\sigma_\zs{Q}$ defined in \eqref{L_1_Q}. 
If it is known, we take $\wh{\sigma}_\zs{n}=\sigma_\zs{Q}$, otherwise, we can choose it, for example, as in \cite{KonevPergamenshchikov2012}, i.e.
\begin{equation}\label{sec:Mo.4+===}
\wh{\sigma}_\zs{n}=
\sum^n_\zs{j=[\sqrt{n}]+1}\,\wh{t}^2_\zs{j,n}\,,
\end{equation}
where $\wh{t}_\zs{j,n}$ are the estimators for the Fourier coefficients with respect to the trigonometric basis \eqref{sec:In.5}, i.e.
\begin{equation*}\label{sec:Mo.4-2-31-3}
\wh{t}_\zs{j,n}=\frac{1}{n}
\int^{n}_\zs{0}\,Tr_\zs{j}(t)\d y_\zs{t}\,.
\end{equation*}

\noindent 
Finally, in order to choose the weights, we will minimize the following cost function
\begin{equation*}\label{sec:Mo.5}
J_n(\lambda)=\sum_\zs{j=1}^{n} \lambda^2(j) \wh{\theta}^2_\zs{j,n} -2 \sum_\zs{j=1}^{n} \lambda(j)\wt{\theta}_\zs{j,n} + \delta\,P_\zs{n}(\lambda),
\end{equation*}
where $\delta>0$ is some  threshold which will be specified later and the penalty term is
\begin{equation}\label{sec:Mo.6}
P_\zs{n}(\lambda)= \frac{  \wh{\sigma}_\zs{n} |\lambda|^2}{n}.
\end{equation}
\noindent
We define the model selection procedure as
\begin{equation}\label{sec:Mo.9-MS}
\wh{S}_* = \wh{S}_\zs{\hat \lambda}
\quad\mbox{with}\quad
\wh{\lambda}= \mbox{argmin}_\zs{\lambda\in\Lambda} J_n(\lambda)\,.
\end{equation}
We recall that the set $\Lambda$ is finite so $\hat \lambda$ exists. In the case when $\hat \lambda$ is not unique, we take one of them.

As is shown in 
\cite{BarbuBeltaifPergamenshchikov2017a, KonevPergamenshchikov2012,  PchelintsevPchelintsev_VPergamenshchikov2017}
 both Conditions $\C_\zs{1})$
and $\C_\zs{2})$ hold for the processes \eqref{sec:In.1+1}, \eqref{sec:Ex.1} and \eqref{sec:Ex.5}.

\begin{proposition}\label{Pr.sec:2.1}
If the conditions $\C_\zs{1})$ and $\C_\zs{2})$ hold for the
distribution $Q$ of the process $\xi$ in \eqref{sec:In.1}, then, for any $n\geq1$ and $0<\delta<1/3$,
the risk \eqref{sec:risks_00} of estimate \eqref{sec:Mo.9-MS} for $S$
satisfies the oracle inequality
\begin{equation}\label{sec:Mo.15_OrIneq}
\cR_\zs{Q}(\wh{S}_\zs{*},S)\,\le\, \frac{1+3\delta}{1-3\delta}
\min_\zs{\lambda\in\Lambda} \cR_\zs{Q}(\wh{S}_\zs{\lambda},S)
+
\frac{\B_\zs{n}(Q)}{\delta n}
\,,
\end{equation}
where 
$\B_\zs{n}(Q)=\U_\zs{n}(Q)
\left( 
1+\vert\Lambda\vert_\zs{*}\E_\zs{Q}|\wh{\sigma}_\zs{n}-\sigma_\zs{Q}|
\right)
$
and
 the coefficient $\U_\zs{n}(Q)$ is such that for any $\epsilon>0$
\begin{equation}
\label{termB_rest}
\lim_\zs{n\to\infty}\frac{\U_\zs{n}(Q)}{n^{\epsilon}}=0
\,.
\end{equation}
\end{proposition}

\noindent
 In the case, when the value of $\sigma_\zs{Q}$ is known, one can take
$\wh{\sigma}_\zs{n}=\sigma_\zs{Q}$ and
\begin{equation*}\label{sec:Mo.9_P}
P_\zs{n}(\lambda)=\frac{\sigma_\zs{Q}\,|\lambda|^2_\zs{n}}{n}\,,
\end{equation*}
then we can rewrite the oracle inequality \eqref{sec:Mo.15_OrIneq}  with $\B_\zs{n}(Q)=\U_\zs{n}(Q)$. Also we study the accuracy properties for the estimator \eqref{sec:Mo.4+===}.

\begin{proposition}\label{sec:Mo.Prop.1}
Let in the model \eqref{sec:In.1} the function $S(\cdot)$ is continuously differentiable.
Then, for any $n\geq 2$,
\begin{equation*}\label{sec:Mo.15_kappa}
\E_\zs{Q}|\wh{\sigma}_\zs{n}-\sigma_\zs{Q}| \leq
\frac{\varkappa_\zs{n}(Q)(1+\|\dot{S}\|^2)}{\sqrt{n}}\,,
\end{equation*}
where the term $\varkappa_\zs{n}(Q)$ possesses  the property  \eqref{termB_rest}.
\end{proposition}
\noindent
To obtain the oracle inequality for the robust risk \eqref{sec:risks}
 we need some additional condition on the distribution family $\cQ_\zs{n}$.
We set
\begin{equation}
\label{sigma*_n}
\varsigma^{*}=\varsigma^{*}_\zs{n}=\sup_\zs{Q\in\cQ_\zs{n}}
\sigma_\zs{Q}
\quad\mbox{and}\quad
\L^{*}_\zs{n}=
\sup_\zs{Q\in\cQ_\zs{n}}(\L_\zs{1,n}(Q)+\L_\zs{2,n}(Q))
\,.
\end{equation}

\noindent $\C^{*}_\zs{1}$) {\it Assume that the conditions $\C_\zs{1})$--$\C_\zs{2})$
hold and for any  $\epsilon>0$
$$
\lim_\zs{n\to\infty}\frac{\L^{*}_\zs{n}+\varsigma^{*}_\zs{n}}{n^{\epsilon}}
=0\,.
$$
}

\noindent 
Now we impose the conditions on the set of the weight coefficients $\Lambda$.

\noindent $\C^{*}_\zs{2})$ {\it Assume that the set $\Lambda$ is such that for any $\epsilon>0$
\begin{equation*}\label{sec:Mo.8+1}
\lim_\zs{n\to\infty}\frac{\nu}{n^{\epsilon}}=0
\quad\mbox{and}\quad
\lim_\zs{n\to\infty}\,\frac{\vert\Lambda\vert_\zs{*}}{n^{1/2+\epsilon}}
=0\,.
\end{equation*}
}

\begin{theorem}\label{Th.sec:2.3}
Assume that the conditions 
$\C^*_\zs{1})$--$\C^*_\zs{2})$ hold.  Then the robust risk
 \eqref{sec:risks} of the estimate \eqref{sec:Mo.9-MS} for
continuously differentiable function $S(t)$ satisfies for any $n\ge 2$  and
$0<\delta<1/3$ the oracle inequality
\begin{align*}\label{sec:Mo.20}
\cR^{*}_\zs{n}(\wh{S}_\zs{*},S)\,\le\,
\frac{1+3\delta}{1-3\delta} \min_\zs{\lambda\in\Lambda}
\cR_n^{*}(\wh{S}_\zs{\lambda},S) +\frac{1}{\delta n}\,\B^{*}_\zs{n}(1+\|\dot{S}\|^2)\,,
\end{align*}
where the term $\B^{*}_\zs{n}$ satisfies the property
\eqref{termB_rest}.
\end{theorem}

Now we specify the weight coefficients $(\lambda(j))_\zs{j\ge
1}$ in the way proposed in \cite{GaltchoukPergamenshchikov2009a}
 for a heteroscedastic regression
model in discrete time. First we define the normalizing coefficient which defined the minimax convergence rate
\begin{equation}
\label{upsilon-nnn}
 v_\zs{n}=
\frac{n}{\varsigma^{*}}
\,,
\end{equation}
where the upper proxy variance is $\varsigma^{*}$ is defined in \eqref{sigma*_n}.  Consider a numerical grid of the form
\begin{equation*}\label{sec:Imp.7}
\cA_\zs{n}=\{1,\ldots,k^*\}\times\{r_1,\ldots,r_m\}\,,
\end{equation*}
where  $r_i=i\varepsilon$ and $m=[1/\varepsilon^2]$. Both
parameters $k^*\ge 1$ and $0<\varepsilon\le 1$ are assumed to be
functions of $n$, i.e. $ k^*=k^*(n)$ and
$\varepsilon=\varepsilon(n)$, such that for any $\delta>0$
\begin{equation*}\label{sec:Imp.8}
\left\{
\begin{array}{ll}
&\lim_\zs{n\to\infty}\,k^*(n)=+\infty\,,
\quad
\lim_\zs{n\to\infty}\,\dfrac{k^*(n)}{\ln n}=0\,,\\[6mm]
&\lim_\zs{n\to\infty}\varepsilon(n)=0
\quad\mbox{and}\quad
\lim_\zs{n\to\infty}\,n^{\delta}\varepsilon(n)\,=+\infty\, .
\end{array}
\right.
\end{equation*}
One can take, for
example,
$$
\varepsilon(n)=\frac{1}{\ln (n+1)}
\quad\mbox{and}\quad
k^*(n)=\sqrt{\ln (n+1)}\,.
$$
For each $\alpha=(\beta,r)\in\cA_\zs{n}$ we introduce the weight
sequence $\lambda_\zs{\alpha}=(\lambda_\zs{\alpha}(j))_\zs{j\ge 1}$
as
\begin{equation}\label{sec:Imp.9}
\lambda_\zs{\alpha}(j)=\Chi_\zs{\{1\le j\le d\}}+
\left(1-(j/\omega_\alpha)^\beta\right)\, \Chi_\zs{\{ d<j\le
\omega_\alpha\}}
\end{equation}
where $d=d(\alpha)=\left[\omega_\zs{\alpha}/\ln (n+1)\right]$,
$
\omega_\zs{\alpha}=\left(\tau_\zs{\beta}\,r\,v_\zs{n}\right)^{1/(2\beta+1)}$ and
$$
\tau_\zs{\beta}=\frac{(\beta+1)(2\beta+1)}{\pi^{2\beta}\beta}\,.
$$
We set
\begin{equation}\label{sec:Imp.10_Lambda}
\Lambda\,=\,\{\lambda_\zs{\alpha}\,,\,\alpha\in\cA_\zs{n}\}\,.
\end{equation}
It will be noted that in this case the cardinal of the set $\Lambda$ is  
$\nu=k^{*} m$. Moreover,
taking into account that $\tau_\zs{\beta}<1$ for $\beta\ge 1$ 
we obtain for the set \eqref{sec:Imp.10_Lambda}
\begin{equation*}
\label{sec:Ga.1++1--2}
 \vert \Lambda\vert_\zs{*}\,
 \le\,1+
\sup_\zs{\alpha\in\cA}
  \omega_\zs{\alpha}
\le 1+(\upsilon_\zs{n}/\ve )^{1/3}\,.
\end{equation*}

Note that the form \eqref{sec:Imp.9}
for the weight coefficients 
 was proposed by Pinsker in \cite{Pinsker1981}
 for the efficient estimation in the nonadaptive case, i.e. when the regularity parameters of the function $S$ are known. 
In the adaptive case  these weight coefficients are
  used in \cite{KonevPergamenshchikov2012, KonevPergamenshchikov2015}
   to show the asymptotic efficiency for model selection procedures.

\section{Improved estimation}\label{sec:Imp}

In this Section we consider the  improved estimation method for the model \eqref{sec:In.1}.
We impose the following additional condition on the noise distribution.

$\D_\zs{1}$) {\sl There exists $n_\zs{0}\ge 1$ such that
for any $n\ge n_\zs{0}$ there exists a  $\sigma$ - field $\cG_\zs{n}$
for which
 the random vector
 $\wt{\xi}_\zs{d,n}=(\xi_\zs{j,n})_\zs{1\le j\le d}$
 is the $\cG_\zs{n}$ conditionally Gaussian in $\bbr^{d}$ with the covariance matrix
\begin{equation*}\label{sec:Imp.6-1}
\G_\zs{n}=\left(
\E\,\xi_\zs{i,n}\,\xi_\zs{j,n}|\cG_\zs{n})
\right)_\zs{1\le i,j\le d}
\end{equation*}
and for some nonrandom
constant $l^{*}_\zs{n}>0$
\begin{equation*}\label{sec:Imp.6-1-0}
\inf_\zs{Q\in\cQ_\zs{n}}\quad
\left(
\tr \,\G_\zs{n}
-
\,
\lambda_\zs{max}(\G_\zs{n})
\right)
\geq l^{*}_\zs{n}\,,
\end{equation*}
where $\lambda_\zs{max}(A)$ is the maximal eigenvalue of the matrix $A$.
}

\noindent

\begin{proposition}\label{sec:Imp.Prop_2_1}
Let in the model \eqref{sec:In.1} the noise process describes by the L\'evy process \eqref{sec:In.1+1}. Then  the
condition $\D_\zs{1})$ holds with $l^{*}_\zs{n}=(d-1)\underline{\varrho}$  for any $n\ge 1$.
\end{proposition}

\begin{proposition}\label{sec:Imp.Prop_2_2}
Let in the model \eqref{sec:In.1} the noise process describes by the Ornstein--Uhlenbeck process \eqref{sec:Ex.1}. Then  the
condition $\D_\zs{1})$ holds with $l^{*}_\zs{n}=(d-6)\underline{\varrho}/2$  for any $n\ge n_0$ and
$d\ge d_0=\inf\{d\ge 7\,:\,5+\ln d\le \check{a} d\}$,
$\check{a}=(1-e^{-a_\zs{max}})/(4a_\zs{max})$.
\end{proposition}

Now,
for the first $d$ Fourier coefficients
 we use
 the improved estimation method proposed
for parametric models
 in
\cite{Pchelintsev2013}. To this end we
set $\wt{\theta}_\zs{n}=(\wh{\theta}_\zs{j,n})_\zs{1\le j\le d}$.
In the sequel we will use the norm $\vert x\vert^{2}_\zs{d}=\sum^{d}_\zs{j=1}\,x^{2}_\zs{j}$
for any vector $x=(x_\zs{j})_\zs{1\le j\le d}$ from $\bbr^{d}$.
Now
we define the shrinkage estimators as
\begin{equation*}\label{sec:Imp.12}
\theta^{*}_\zs{j,n}=
\left(1-g(j)\right)\wh{\theta}_\zs{j,n}\,,
\end{equation*}
where  $g(j)=(\c_\zs{n}/|\wt{\theta}_\zs{n}|_\zs{d}) \Chi_\zs{\{1\le j\le d\}}$ and
\begin{equation}
\label{coefficient_c_n}
\c_\zs{n}=
\frac{l^{*}_\zs{n}}{\left(r^{*}_\zs{n}+\sqrt{d/v_\zs{n}}\right)\,n}\,.
\end{equation}
The positive parameter
 $r^{*}_\zs{n}$  is such that $\lim_\zs{n\to\infty}\,r^{*}_\zs{n}\,=\infty$ and
for any $\epsilon>0$
\begin{equation}\label{sec:Imp.12+1_r_n}
\lim_\zs{n\to\infty}\,
\frac{r^{*}_\zs{n}}{n^{\epsilon}}
\,=\,0
\end{equation}
and $v_\zs{n}$ defined in \eqref{upsilon-nnn}. Now we introduce a class of shrinkage
weighted least squares estimates for $S$ as
\begin{equation}\label{sec:Imp.11}
S^{*}_\zs{\lambda}=\sum^{n}_\zs{j=1}\lambda(j)\theta^{*}_\zs{j,n}\phi_\zs{j}\,.
\end{equation}

We denote the difference of quadratic risks of the estimates \eqref{sec:Mo.1+??+} and \eqref{sec:Imp.11} as
$$
\Delta_{Q}(S):=\cR_\zs{Q}(S^{*}_\zs{\lambda},S)-\cR_\zs{Q}(\wh{S}_\zs{\lambda},S)\,.
$$
We obtain the following result.

\begin{theorem}\label{Th.sec:Imp.1}
Let the observed process $(y_\zs{t})_\zs{0\le t \le n}$ describes by the equation \eqref{sec:In.1} and the condition
$\D_\zs{1})$ holds. Then for any $n\ge 1$
\begin{equation}\label{sec:Imp.11+1}
\sup_{Q\in\cQ_\zs{n}}\,\sup_\zs{\Vert S\Vert\le r^{*}_\zs{n}}
\Delta_{Q}(S)\le-\c^2_\zs{n}
\,.
\end{equation}
\end{theorem}

\begin{remark}\label{Re;sec:Imp.1}
The inequality \eqref{sec:Imp.11+1} means that non asymptotically, i.e. for any $n\ge 1$,
 the estimate \eqref{sec:Imp.11}   outperforms in mean square accuracy the estimate \eqref{sec:Mo.1+??+}.
Moreover in the efficient weight coefficients  $d \approx n^{\check{\delta}}$ as $n\to \infty$ for some $\check{\delta}>0$.
Therefore, in view of the  definition  \eqref{coefficient_c_n} and the conditions
\eqref{sec:Ex.01-2}
 and
 \eqref{sec:Imp.12+1_r_n}
 $n \c_\zs{n}\to\infty$ as $n\to\infty$.
This means that improvement is considerably may better than for the parametric regression
when the parameter dimension $d$ is fixed  \cite{Pchelintsev2013}.
\end{remark}

\bigskip

\section{Improved model selection}\label{sec:IMo}

This Section gives the construction of a model selection procedure  for
estimating a function $S$ in \eqref{sec:In.1} on the basis of improved weighted least square estimates $(S^{*}_\zs{\lambda})_\zs{\lambda\in\Lambda}$ and states
the sharp oracle inequality for the robust risk of proposed procedure.

As in Section 3, the performance of any estimate $S^{*}_\zs{\lambda}$ will be measured by the
empirical squared error
$$
\Er_\zs{n}(\lambda)=\|S^*_\zs{\lambda}-S\|^2.
$$
In order to obtain a good estimate, we have to write a rule to choose a weight vector
$\lambda\in \Lambda$ in \eqref{sec:Imp.11}. It is obvious, that the best way is to minimise
the empirical squared error with respect to $\lambda$. Making use the estimate definition
\eqref{sec:Imp.11} and the Fourier transformation of $S$ implies
\begin{equation*}\label{sec:Mo.1+11+}
\Er_\zs{n}(\lambda)\,=\,
\sum^{n}_\zs{j=1}\,\lambda^2(j)(\theta^*_\zs{j,n})^2\,-
2\,\sum^{n}_\zs{j=1}\,\lambda(j)\theta^*_\zs{j,n}\,\theta_\zs{j}\,+\,
\sum^{n}_\zs{j=1}\theta^2_\zs{j}\,.
\end{equation*}
Here one needs to replace  the terms
$\theta^*_\zs{j,n}\,\theta_\zs{j}$ by their estimators
$\bar{\theta}_\zs{j,n}$. We set
\begin{equation*}\label{sec:Mo.2+=???}
\bar{\theta}_\zs{j,n}=
\theta^*_\zs{j,n}\,\wh{\theta}_\zs{j,n}-\frac{\wh{\sigma}_\zs{n}}{n}\,,
\end{equation*}
where $\wh{\sigma}_\zs{n}$ is defined in \eqref{sec:Mo.4+===}.
For this change in the empirical squared error, one has to pay
some penalty. Thus, one comes to the cost function of the form
\begin{equation*}\label{sec:Mo.4+++???}
J_\zs{n}^*(\lambda)\,=\,\sum^{n}_\zs{j=1}\,\lambda^2(j)(\theta^*_\zs{j,n})^2\,-
2\,\sum^{n}_\zs{j=1}\,\lambda(j)\,\bar{\theta}_\zs{j,n}\,
+\,\delta\,P_\zs{n}(\lambda)\,,
\end{equation*}
where $\delta$ is some positive constant and the penalty term
$P_\zs{n}(\lambda)$ is  defined  in \eqref{sec:Mo.6}.
Substituting the weight coefficients, minimizing the cost function
\begin{equation}\label{sec:Mo.6-1}
\lambda^*=\mbox{argmin}_\zs{\lambda\in\Lambda}\,J_n^*(\lambda)\,,
\end{equation}
in \eqref{sec:Imp.11} leads to the improved model selection procedure
\begin{equation}\label{sec:Mo.7}
S^*=S^*_\zs{\lambda^*}\,.
\end{equation}
It will be noted that $\lambda^*$ exists because
 $\Lambda$ is a finite set and also if the
minimizing sequence in \eqref{sec:Mo.6-1} $\lambda^*$ is not
unique, one can take any minimizer.

\begin{theorem}\label{sec:Mo.Th.1}
If the conditions $\C_\zs{1})$ and $\C_\zs{2})$ hold for the
distribution $Q$ of the process $\xi$ in \eqref{sec:In.1}, then, for any $n\geq1$ and $0<\delta<1/3$,
the risk \eqref{sec:risks_00} of estimate \eqref{sec:Mo.7} for $S$
satisfies the oracle inequality
\begin{equation*}\label{sec:Mo.10}
\cR_\zs{Q}(S^*_\zs{\lambda^*},S)\,\le\, \frac{1+3\delta}{1-3\delta}
\min_\zs{\lambda\in\Lambda} \cR_\zs{Q}(S^*_\zs{\lambda},S)
+\frac{\check{\B}_\zs{n}(Q)}{n \delta}\,,
\end{equation*}
where 
$\check{\B}_\zs{n}(Q)=\check{\U}_\zs{n}(Q)
\left( 
1+\vert\Lambda\vert_\zs{*}\E_\zs{Q}|\wh{\sigma}_\zs{n}-\sigma_\zs{Q}|
\right)
$
and
 the coefficient $\check{\U}_\zs{n}(Q)$ satisfies the property
\eqref{termB_rest}.
\end{theorem}

\vspace{2mm}

Now Theorem \ref{sec:Mo.Th.1} and Proposition \ref{sec:Mo.Prop.1} directly
imply the following inequality for the robust risk \eqref{sec:risks} of the procedure  \eqref{sec:Mo.7}.

\begin{theorem}\label{sec:Mo.Th.2}
Assume that the conditions $\C_\zs{1}^*)$ and
$\C_\zs{2}^*)$ hold and the function $S$ is
continuously differentiable.
 Then  for any $n\ge 2$  and
$0<\delta<1/3$ 
\begin{equation*}\label{sec:Mo.10+1-2}
\cR^{*}_\zs{n}
(S^*_\zs{\lambda^*},S)\,\le\, \frac{1+3\delta}{1-3\delta}
\min_\zs{\lambda\in\Lambda} \cR^{*}_\zs{n}\,(S^*_\zs{\lambda},S)
+\frac{\check{\U}^{*}_\zs{n}(1+\|\dot{S}\|^2)}{n \delta}\,,
\end{equation*}
where 
 the coefficient $\check{\U}^{*}_\zs{n}$ satisfies the property
\eqref{termB_rest}.
\end{theorem}

\section{Asymptotic efficiency}\label{sec:Ae}

In order to study the asymptotic efficiency we define the following functional Sobolev ball
\begin{equation*}\label{sec:Ae.1}
W_\zs{k,\r}=\{f\in\C^{k}_\zs{p}[0,1]\,:\,
\sum_\zs{j=0}^k\,\|f^{(j)}\|^2\le \r\}\,,
 \end{equation*}
where $\r>0$ and $k\ge 1$ are
some unknown parameters, $\C^{k}_\zs{p}[0,1]$ is the space of
 $k$ times differentiable $1$ - periodic $\bbr\to\bbr$ functions
 such that $f^{(i)}(0)=f^{(i)}(1)$  for any $0\le i \le k-1$. It is well known that for any $S\in W_\zs{k,\r}$
 the optimal rate of convergence is
$n^{-2k/(2k+1)}$ (see, for example, \cite{Pinsker1981, Nussbaum1985}).
On the basis of the model selection procedure
we construct the adaptive
procedure $S_\zs{*}$ for which we obtain the following asymptotic upper bound
for the quadratic risk, i.e.
 we show that the parameter \eqref{sec:Ae.3} gives a lower bound
for the asymptotic normalized risks.
To this end we denote by $\Sigma_\zs{n}$ of all estimators $\wh{S}_\zs{n}$ of $S$ measurable with respect to
the process
\eqref{sec:In.1}, i.e.
$\sigma\{y_\zs{t}\,,\,0\le t\le n\}$.

\begin{theorem}\label{Th.sec: Ae.2}
The robust risk \eqref{sec:risks}
admits the following asymptotic lower bound
 \begin{equation*}\label{sec:Ae.5++==+??}
\liminf_\zs{n\to\infty}\,
\inf_\zs{\wh{S}_\zs{n}\in\Sigma_\zs{n}}
\,v_\zs{n}^{2k/(2k+1)}
\sup_\zs{S\in W_\zs{k,r}}\,\cR^{*}_\zs{n}(\wh{S}_\zs{n},S)
\,
\ge l_\zs{*}(\r)  \,.
 \end{equation*}
\end{theorem}

\noindent
We show that this lower bound is sharp in the following sense.

\begin{theorem}\label{Th.sec: Ae.1}
The quadratic risk \eqref{sec:risks} for the
 estimating procedure $S^{*}$ has the following asymptotic upper bound
 \begin{equation*}\label{sec:Ae.4}
\limsup_\zs{n\to\infty}\,v_\zs{n}^{2k/(2k+1)}
\sup_\zs{S\in W_\zs{k,\r}}\,\cR^{*}_\zs{n}(S^{*},S)
\,
\le l_\zs{*}(\r)
\,.
 \end{equation*}
\end{theorem}

\noindent
It is clear that Theorem \ref{Th.sec: Ae.1}
 and Theorem \ref{Th.sec: Ae.2}
imply
\begin{corollary}\label{Co.sec: Ae.1}
The model selection procedure $S^{*}$
is efficient, i.e.
\begin{equation}\label{sec:Ae.5+6+?}
\lim_{n\to\infty}\,(v_\zs{n})^{\frac{2k}{2k+1}}\,
\sup_\zs{S\in W_\zs{k,\r}}\,\cR^{*}_\zs{n}(S^{*},S)\,
= l_\zs{*}(\r)
\,.
\end{equation}
\end{corollary}

\noindent 
Note that the equality 
\eqref{sec:Ae.5+6+?} implies that the parameter
\eqref{sec:Ae.3} is the Pinsker constant
in this case (cf. \cite{Pinsker1981}).
Moreover, it  means that the robust efficiency holds with
the  convergence rate $(v_\zs{n})^{\frac{2k}{2k+1}}$.
It is well known that for the simple risks
 the optimal (minimax) estimation convergence rate
for the functions from the set  $W_\zs{k,r}$
 is $n^{2k/(2k+1)}$ (see, for example, \cite{IbragimovKhasminskii1981, Nussbaum1985, Pinsker1981}).
So,
 if  the distribution upper bound  $\varsigma^{*}\to 0$  as $n\to\infty$
 we obtain the more rapid rate, and
 if $\varsigma^{*}\to \infty$  as $n\to\infty$
we obtain the more slow rate. In the case when $\varsigma^{*}$ is constant the robust rate is the same as the classical non robust convergence rate.

\bigskip

\noindent 
{\bf Acknowledgements.}

The first author is partially supported by the 
 the RSF grant 17-11-01049.
The last author is partially supported  by the Russian Federal Professor program (project no. 1.472.2016/1.4,
  Ministry of Education and Science) and
  by the project XterM - Feder, University of Rouen, France.

\end{document}